\documentclass[a4paper,12pt]{article}
\usepackage{amsmath,amsfonts,amssymb}
\usepackage[all]{xy}
\newtheorem{theorem}{Theorem}[section]
\newtheorem{fact}[theorem]{Fact}
\newtheorem{lemma}[theorem]{Lemma}
\newtheorem{corollary}[theorem]{Corollary}
\newcommand{\Hom}{\textnormal{Hom}}
\newcommand{\Ext}{\textnormal{Ext}}

\title{Another point in homological algebra: 
Duality for discontinuous group actions}
\author{Fritz Grunewald and Wilhelm Singhof\\
{\small Heinrich-Heine-Universit\"at D\"usseldorf}\\
{\small e-mail: fritz@math.uni-duesseldorf.de}
}

\begin{document}

\maketitle

\hrule
\begin{abstract}
We consider discontinuous operations of a group $G$ 
on a contractible $n$-dimensional manifold $X$.
Let $E$ be a finite dimensional representation of $G$ 
over a field $k$ of characteristics $0$.
Let $\mathcal{E}$ be the sheaf on the quotient space 
$Y=G \setminus X$ associated to $E$. 
Let $H^{\bullet}_{\textbf{!}}(Y;\mathcal{E})$ 
be the image in  $H^{\bullet}(Y;\mathcal{E})$
of the cohomology with compact support. In the cases where
both  $H^{\bullet}_{\textbf{!}}(Y;\mathcal{E})$ and 
 $H^{\bullet}_{\textbf{!}}(Y;\mathcal{E}^*)$ ($\mathcal{E}^*$ being the 
the sheaf associated to the representation dual to $E$)
are finite dimensional,
we establish a non-degenerate duality between 
$H^{m}_{\textbf{!}}(Y;\mathcal{E})$  and
$H^{n-m}_{\textbf{!}}(Y;\mathcal{E}^{\ast})$.
We also show that this duality is compatible with Hecke operators.

\end{abstract}
\medskip
\hrule

\tableofcontents

\section{Introduction}

We consider discontinuous operations of a group $G$ 
on a contractible $n$-dimensional manifold $X$.
Let $E$ be a finite dimensional representation of $G$ 
over a field $k$ of characteristics $0$, say.
Associated to $E$ is a sheaf $\mathcal{E}$ 
on the quotient space $Y=G \setminus X$, and a theorem
of Grothendieck \cite{G} says that there is an isomorphism
\begin{equation}
 H^{\bullet}(G;E) \cong H^{\bullet}(Y;\mathcal{E})\; .
\end{equation}
This is a result which is basic for the cohomology 
theory of groups. In $H^{\bullet}(Y;\mathcal{E})$, we have the
canonical subgroup $H^{\bullet}_{\textbf{!}}(Y;\mathcal{E})$ 
which is the image of the cohomology with compact supports.
In Section 2, we show that there is a non-degenerate pairing  
\begin{equation}
B:\ H^{m}_{\textbf{!}}(Y;\mathcal{E})\otimes H^{n-m}_{\textbf{!}}(Y;\mathcal{E}^{\ast})\rightarrow k \; .
\end{equation}
This of course gives rise to isomorphisms 
\begin{equation}\label{duaa}
H^{m}_{\textbf{!}}(Y;\mathcal{E})^* \cong 
H^{n-m}_{\textbf{!}}(Y;\mathcal{E}^{\ast})
\end{equation}
where $\mathcal{E}^{\ast}$ is the sheaf associated 
to the dual representation $E^{\ast}$ and where 
$H^{m}_{\textbf{!}}(Y;\mathcal{E})^*$ is the dual of 
$H^{m}_{\textbf{!}}(Y;\mathcal{E})$.

Special instances of this form of Poincar\'e duality 
were frequently used. In full generality, we derive it
by using Grothendieck's arguments together with Poincar\'e-Verdier duality.
In fact we feel that our paper is another (perhaps missing)
chapter in the work \cite{G} of Grothendieck.
\par
In Section 3, we show that 
the above duality isomorphism (\ref{duaa}) 
is compatible with Hecke operators.

\section{Poincar\'e duality}

Let $X$ be a connected, orientable, $n$-dimensional manifold and let $G$ be a discrete group operating
effectively and smoothly on $X$ preserving the orientation. We assume that this operation is \emph{discontinuous},
i.e. the following two conditions are satisfied:
\begin{itemize}
 \item  For each $x \in X$, the stabilizer $G_{x}$ is finite.
 \item  Each $x$ admits a neighborhood $V$ such that 
\[\gamma \cdot V \cap V = \emptyset \quad \text{for all} \quad \gamma \in G \smallsetminus G_{x}\;.\]
\end{itemize}
 
We write $Y:=G\setminus X$ and denote by $p:X\rightarrow Y$ the natural projection. Let $X_{\textnormal{reg}}$
be the open subset of $X$ consisting of all points with trivial stabilizer.\par
Since the smooth operation of a finite group looks, in the neighborhood of a fixed point of the group, like a
linear operation [B1], we have the following two properties:
\begin{fact}
$X \smallsetminus X_{\textnormal{reg}}$ has at least codimension $2$ in $X$.
\end{fact}

\begin{fact}
 Every point of $Y$ has an open neighborhood $U$ such that there exists an open subset $V$ of $X$ diffeomorphic
to $\mathbb{R}^{n}$ and a finite subgroup $H$ of $G$ such that $H \cdot V=V$ and that $p$ induces a homeomorphism
from $H \setminus V$ onto $U$.
\end{fact}

Since paths in quotients by finite groups can be lifted ([B1], Th. 6.2, p. 91), we conclude from 2.2:
\begin{fact}
 Let $A$ be a pathwise connected subset of $Y$ and let $A_{0}$ be a path component of $p^{-1}(A)$. Then:
\begin{enumerate}
 \item [a)] $p^{-1}(A)=G \cdot A_{0}$.
\item[b)] For each $\gamma \in G$, either $\gamma \cdot A_{0}=A_{0}$ or $\gamma \cdot A_{0}\cap A_{0}=\emptyset$. 
\end{enumerate}

\end{fact}

Let $k$ be a field with the following property: For each $x \in X$, the characteristics of $k$ is prime to
the order of $G_{x}$. \par
With the terminology of [B2], Def. V.9.1 and using [B2], Th. II.19.2 and [B2], Th. V.19.3, we then have:
\begin{fact}
 $Y$ is an orientable $n$-dimensional homology manifold over $k$.
\end{fact}

By a representation of $G$, we mean a finite dimensional vector space $E$ over $k$ with a linear operation
of $G$ on it. A sheaf will always be a sheaf of $k$-vector spaces.\par
We consider the following functor $\mathcal{S}$ from the category of representations of $G$ to the
category of sheaves on $Y$: For a representation $E$, the sheaf $\mathcal{E}:=\mathcal{S}(E)$ is given
by 
\[\Gamma (U;\mathcal{E}):= \{\text{locally constant} \; G \text{-equivariant maps} \; p^{-1}(U)\rightarrow E \}\]
for open subsets $U$ of $Y$. By 2.3, the stalks of $\mathcal{E}$ are given by 
\[ \mathcal{E}_{y} =\{G \text{-equivariant maps} \; p^{-1}(y)\rightarrow E \}\, .\]

The constant sheaf on $X$ with stalks $E$, denoted by $E_{X}$, is a special instance of a $G$-sheaf on $X$.
For any $G$-sheaf $\mathcal{F}$ on $X$, the direct image sheaf $p_{\ast}(\mathcal{F})$ is a $G$-sheaf on $Y$,
where $G$ operates trivially on $Y$. Let $p_{\ast}^{G}(\mathcal{F})$ be the subsheaf of $p_{\ast}(\mathcal{F})$
consisting of those elements which are fixed under $G$. With this notation, we have
\[\mathcal{S}(E)=p_{\ast}^{G}(E_{X})\; .\]
More generally, following [G] (5.2.4), we introduce the sheaves 
\[\mathcal{H}^{m}(G;\mathcal{F}):=R^{m}p_{\ast}^{G}(\mathcal{F})\]
on $Y$ so that $\mathcal{S}(E)=\mathcal{H}^{0}(G;E_{X})$.

\medskip
We denote by $\Gamma_{X}^{G}(\mathcal{F})$ the $G$-invariant sections of $\mathcal{F}$ and put as in [G] (5.2.3):
\[ H^{m}(X;G,\mathcal{F}):=R^{m}\Gamma_{X}^{G}(\mathcal{F})\; .\]
Moreover, let $\Gamma_{X,c}^{G}(\mathcal{F})$ denote the elements in $Gamma_{X}^{G}(\mathcal{F})$ with
compact support and put as in [G] (5.7.2):
\[H^{m}_{c}(X;G,\mathcal{F}):= R^{m}\Gamma_{X,c}^{G}(\mathcal{F})\; .\]
By [G], Th. 5.2.1 there are two spectral sequences

\begin{equation}
 (H^{p}(Y;\mathcal{H}^{q}(G;\mathcal{F})))_{p,q}\Rightarrow H^{\bullet}(X;G,\mathcal{F})\; ,
\end{equation}
\begin{equation}
 (H^{p}(G;H^{q}(X;\mathcal{F})))_{p,q}\Rightarrow H^{\bullet}(X;G,\mathcal{F})\; .
\end{equation}

\medskip
For a \emph{finite} group $G$, there are, by [G] (5.7.4) and (5.7.7), analogs for the cohomology with
compact supports.

\begin{equation}
 (H^{p}_{c}(Y;\mathcal{H}^{q}(G;\mathcal{F})))_{p,q}\Rightarrow H^{\bullet}_{c}(X;G,\mathcal{F})\; ,
\end{equation}
\begin{equation}
 (H^{p}(G;H^{q}_{c}(X;\mathcal{F})))_{p,q}\Rightarrow H^{\bullet}_{c}(X;G,\mathcal{F})\; .
\end{equation}

\medskip
From the spectral sequences (3) and (4), Grothendieck deduces the basic relation (1). Similarly, from
(5) and (6) we get:
\begin{fact}
 If the group $G$ is finite and if \textnormal{char} $k$ is prime to the order of $G$, we have
\[H_{c}^{m}(Y;\mathcal{E})\cong H^{m}_{c}(X;E)^{G}\; .\]
\end{fact}

\medskip

For two sheaves $\mathcal{F,G}$ on $Y$ and $p\geq 0$, we have the sheaf $\mathcal{E}xt^{p}(\mathcal{F,G})$.
It is the sheafification of the presheaf
\[U\mapsto \Ext^{p}(\mathcal{F}| U,\mathcal{G}| U)\; . \]

We write $\mathcal{H}om$ instead of $\mathcal{E}xt^{0}$. Observe thar the presheaf
\[U\mapsto \Hom (\mathcal{F}| U,\mathcal{G}| U) \]
is a sheaf, hence $\mathcal{H}om(\mathcal{F,G})(U)=\Hom (\mathcal{F}| U,\mathcal{G}| U)$.

\begin{lemma}
 For two representations $E,F$ of $G$ and their associated sheaves $\mathcal{E,F}$, there is a canonical isomorphism
\[\mathcal{S}(\Hom_{k}(E,F))\cong \mathcal{H}om(\mathcal{E,F})\, .\]
\end{lemma}

\emph{Proof. } Writing $\mathcal{G}:=\mathcal{S}(\Hom_{k}(E,F))$, we have a canonical map
\[\alpha :\mathcal{G} \rightarrow \mathcal{H}om(\mathcal{E,F})\; .\]
The non-trivial point is to show that
\[\alpha_{y} :\mathcal{G}_{y} \rightarrow \mathcal{H}om(\mathcal{E,F})_{y}\]
is surjective for all $y \in Y$. We have 
\[\mathcal{G}_{y}=\{ G \text{-equivariant maps}\; \varphi :p^{-1}(y)\rightarrow \Hom_{k}(E,F)\}\; .\]
Fixing a point $x \in p^{-1}(y)$, the evaluation $\varphi \mapsto \varphi_{x}$ gives an isomorphism
\[\mathcal{G}_{y}\cong \Hom_{kG_{x}}(E,F)\; .\]
Given an element $\delta \in \mathcal{H}om(\mathcal{E,F})_{y}$, we choose an open neighborhood
$U$ of $y$ in $Y$ which is so small that each component of $p^{-1}(U)$ contains exactly one point of 
$p^{-1}(y)$ and that $\delta$ can be represented by an element $\psi \in \Hom (\mathcal{E}|U,\mathcal{F}|U)$.
Let $U_{0}$ be the component of $p^{-1}(U)$ containing $x$. Let $\xi \in U_{0}\cap X_{\text{reg}}$ and $\eta :=p(\xi)$.
Evaluation in $\xi$ gives isomorphisms $\mathcal{E}_{\eta} \cong E$ and $\mathcal{F}_{\eta} \cong F$. Hence the
homomorphism $\psi_{\eta}: \mathcal{E}_{\eta}\rightarrow \mathcal{F}_{\eta}$ gives a $k$-homomorphism
$\varphi :E\rightarrow F$. Using 2.1, it is easy to see that $\varphi$ doesn't depend on the choice of $\xi$.
Moreover, this implies that $\varphi \in \Hom_{kG_{x}}(E,F) \cong \mathcal{G}_{y}$. It is clear that 
$\alpha_{y}(\varphi )=\delta$. \hfill $\Box$

\medskip 

Given a representation $E$ of $G$, we denote by $E^{\ast}:=\Hom_{k}(E,k)$ the dual representation and write
$\mathcal{E}^{\ast}:=\mathcal{S}(E^{\ast})$.

\begin{corollary}
 $\quad\mathcal{H}om(\mathcal{E,F}) \cong \mathcal{H}om(\mathcal{F}^{\ast},\mathcal{E}^{\ast})$.
\end{corollary}

Because of 2.4, we can apply Poincar\'e-Verdier duality 
(see e.g. [GM], (III. 70)) and obtain:
\begin{equation}
 \Hom_{k}(H^{m}_{c}(Y;\mathcal{E}),k) \cong \Ext^{n-m}(\mathcal{E},k_{Y}).
\end{equation}

\begin{lemma}
 $\quad\mathcal{E}xt^{p}(\mathcal{E},k_{Y})=0$ for $p>0$.
\end{lemma}

\emph{Proof. } This is a local property. By 2.2 it suffices to assume that $X$ is diffeomorphic to $\mathbb{R}^{n}$
and that $G$ is finite. We have 
\[\begin{array}{rcll}
 \Ext^{p}(\mathcal{E},k_{Y})& \cong & \Hom_{k}(H_{c}^{n-p}(Y;\mathcal{E}),k) & \text{by (7)}\\
                            & \cong & \Hom_{k}(H_{c}^{n-p}(X;E)^{G},k) & \text{by 2.5.}
\end{array}\]
Since $X\approx\mathbb{R}^{n}$, we have $H_{c}^{n-p}(X;E)=0$ for $p>0$. {\hfill}$\Box$

\begin{lemma}
 For each $p\geq 0$, there is a canonical isomorphism 
\[\Ext^{p}(\mathcal{E},k_{Y}) \cong H^{p}(Y;\mathcal{E}^{\ast})\; .\]
\end{lemma}

\emph{Proof. }  We use the spectral sequence ([G], Th. 4.2.1)
\[(H^{p}(Y;\mathcal{E}xt^{q}(\mathcal{E},k_{Y})))_{p,q}\Rightarrow \Ext^{\bullet}(\mathcal{E},k_{Y})\; .\]
By 2.8, we obtain $\Ext^{p}(\mathcal{E},k_{Y}) \cong H^{p}(Y;\mathcal{H}om(\mathcal{E},k_{Y}))$.
Then apply 2.7. \hfill $\Box$

\medskip
There is a canonical homomorphism of sheaves
\[\mu :\mathcal{E} \otimes \mathcal{E}^{\ast}\rightarrow k_{Y}\]
which induces 
\[\mu_{\ast}:H^{n}_{c}(Y;\mathcal{E} \otimes \mathcal{E}^{\ast})\rightarrow H^{n}_{c}(Y;k)\; .\]
Composing with the cup-product
\[\cup: H^{m}_{c}(Y;\mathcal{E})\otimes H^{n-m}(Y;\mathcal{E}^{\ast})\rightarrow H^{n}_{c}(Y;\mathcal{E} \otimes \mathcal{E}^{\ast})\]
and with
\[\int:H_{c}^{n}(Y;k)\rightarrow k\; ,\]
we obtain for each $m\geq 0$ a homomorphism
\[B:H^{m}_{c}(Y;\mathcal{E})\otimes H^{n-m}(Y;\mathcal{E}^{\ast}) \rightarrow k\, .\]
Combining Poincar\'e-Verdier duality (7) and 
in particular the remark in [GM] following the formula (III.70)
with 2.9, we obtain 
\begin{lemma}
 The pairing $B$ is non-singular for each $m\geq 0$. \hfill $\Box$
\end{lemma}

For a sheaf $\mathcal{F}$ on $Y$, let us write
\[H^{\bullet}_{\textbf{!}}(Y;\mathcal{F}):=\text{im}(H^{\bullet}_{c}(Y;\mathcal{F}) \rightarrow H^{\bullet}(Y;\mathcal{F}))\; .\]

The properties of the cup-product, in particular its anti-commutativity, show that 2.10 implies the following
duality for $H^{\bullet}_{\textbf{!}}$:
\begin{theorem}
 Let $X$ be an orientable, connected, $n$-dimensional manifold and let $G$ be a group operating smoothly
and discontinuously on $X$ preserving the orientation. Let $k$ be a field such that for each $x \in X$, the
characteristcs of $k$ is prime to the order of the stabilizer $G_{x}$. Let $E$ be a representation of $G$
on a finite-dimensional $k$-vector space and let $E^{\ast}$ be the dual representation. We denote by
$\mathcal{E}$ the sheaf on $Y=G\setminus X$ associated to $E$ and by $\mathcal{E}^{\ast}$ the sheaf associated
to $E^{\ast}$ and assume that all cohomology groups $H^{m}(Y;\mathcal{E})$ and $H^{m}(Y;\mathcal{E}^{\ast})$
are finite dimensional $k$-vector spaces. Then $B$ induces for each $m\geq 0$ a non-singular pairing of
finite dimensional $k$-vector spaces
\begin{equation}
B:\ H^{m}_{\textbf{!}}(Y;\mathcal{E})\otimes H^{n-m}_{\textbf{!}}(Y;\mathcal{E}^{\ast})\rightarrow k \; .
\end{equation}
\end{theorem}

\section{Hecke operators}

Let us now assume in addition that there is a group $\mathfrak{G}$ with
\[G \subseteq \mathfrak{G} \subseteq \text{Diff} \,X\]
and that $E$ is actually a representation of $\mathfrak{G}$. Let $g \in \mathfrak{G}$ be an element in
the commensurator of $G$, that is, we assume that the two groups 
\[\begin{array}{rcl}
 \Delta' &:=& gGg^{-1} \cap G,\\
 \Delta'' &:=& G \cap g^{-1}Gg=g^{-1} \Delta'g
\end{array}\]
are of finite index in $G$. Consider the following diagram:

\medskip
\parbox{12cm}{
\begin{equation*}
 \xymatrix{X \ar[rr]^{\displaystyle p'} \ar[rrdd]_{\displaystyle p}&&\Delta'\setminus X \ar[rr]^{\displaystyle\varphi}\ar[dd]^{\displaystyle\pi'}&&\Delta''\setminus X \ar[dd]^{\displaystyle\pi''}
&& X\ar[ll]_{\displaystyle p''}\ar[lldd]^{\displaystyle p}\\
\\
&& G\setminus X && G\setminus X }
\end{equation*}}\hfill
\parbox{8mm}{\begin{equation} \end{equation}}
All the arrows except $\varphi$ are natural projections; $\varphi=\varphi_{g}$ is the diffeomorphism
given by
\[\varphi(\Delta'x):=\Delta''g^{-1}x=g^{-1}\Delta'x \; .\]
As we shall explain, we obtain from (8) homomorphisms

\medskip
\parbox{12cm}{
\begin{equation*}
\xymatrix{H^{\bullet}(\Delta'\setminus X;\mathcal{E}') \ar[dd]^{\displaystyle\tau}&& H^{\bullet}(\Delta''\setminus X;\mathcal{E}'')
\ar[ll]_{\displaystyle \varphi^{\ast}}\\
\\
H^{\bullet}(G\setminus X;\mathcal{E})&& H^{\bullet}(G\setminus X;\mathcal{E}) \ar[uu]_{\displaystyle \pi''^{\ast}}} 
\end{equation*}}\hfill
\parbox{8mm}{\begin{equation} \end{equation}}
Here $\mathcal{E}$ is the sheaf introduced in section 2, and $\mathcal{E}'$ (resp. $\mathcal{E}''$) are the 
analoguous sheaves on $\Delta' \setminus X$ (resp. $\Delta'' \setminus X$).

\textbf{Definition. } The \emph{Hecke operator} $T(g)$ is the endomorphism of \linebreak $H^{\bullet}(G\setminus X; \mathcal{E})$
given by 
\[T(g):= \tau \circ \varphi^{\ast}\circ \pi''^{\ast}\; .\]

Let us now explain the three homomorphisms $\tau, \; \varphi^{\ast}$, and $\pi''^{\ast}$ in (9):\\
We have a canonical inclusion $i':\mathcal{E}\hookrightarrow \pi'_{\ast}\mathcal{E}^{\ast}$. We obtain induced
homomorphisms 
\[\pi'^{\dagger}:H^{\bullet}(G\setminus X; \pi'_{\ast}\mathcal{E}')\rightarrow H^{\bullet}(\Delta'\setminus X; \mathcal{E}')\]
which is an isomorphism by the Vietoris mapping theorem ([B2], Th. II. 11.1) and 
\[\pi'^{\ast}=\pi'^{\dagger}\circ i'_{\ast}:H^{\bullet}(G\setminus X; \mathcal{E})\rightarrow H^{\bullet}(\Delta'\setminus X;
 \mathcal{E}')\; .\]
Similarly for $\mathcal{E}''$ instead of $\mathcal{E}'$. This explains $\pi''^{\ast}$.

\medskip
Next, there is a sheaf homomorphism 
\[\sigma: \pi'_{\ast}\mathcal{E}'\rightarrow \mathcal{E} \]
defined as follows: Choose elements $\gamma_{1}, \ldots , \gamma_{m} \in G$ with 
\[G= \coprod_{i=1}^{m}\gamma_{i}\Delta'\; .\]
For $U$ open in $G\setminus X$, define
\[\sigma_{U}:\Gamma(U;\pi'_{\ast}\mathcal{E}')\rightarrow \Gamma(U;\mathcal{E})\]
\[\sigma_{U}(f)(x):=\sum_{i=1}^{m}\gamma_{i}\cdot f(\gamma_{i}^{-1}\cdot x)\]
for $x \in p^{-1}(U)$. The definition of $\sigma$ doesn't depend on the choice of the $\gamma_{i}$. We then
have a \emph{transfer}
\begin{equation}
 \tau:=\sigma_{\ast}\circ(\pi'^{\dagger})^{-1}:H^{\bullet}(\Delta'\setminus X;\mathcal{E}')\rightarrow H^{\bullet}(G\setminus X;
 \mathcal{E})\; .
\end{equation}
Observe that $\tau$ is a special case of the transfer introduced in [B2], p.140-141. (Bredon's assumption that $G$ is
a finite group is not at all used for the definition of the transfer called $\mu^{\ast}_{G/H}$ by him.)

\medskip
Finally, to explain $\varphi^{\ast}$, we have to specify a $\varphi$-cohomomorphism
\[\Phi:\mathcal{E}''\rightsquigarrow \mathcal{E}'\]
(using the terminology of [B2], I. 4): For each $z \in \Delta'\setminus X$, we need a homomorphism 
\[\Phi_{z}:\mathcal{E}''_{\varphi(z)}\rightarrow \mathcal{E}'_{z}\; .\] 
Choosing an element $x \in p'^{-1}(z)$, we have
\[\mathcal{E}''_{\varphi(z)}=\{ \Delta''\;\text{-equivariant maps}\; f:g^{-1}\Delta'x\rightarrow E\}\; ,\]
\[\mathcal{E}'_{z}=\{ \Delta'\;\text{-equivariant maps}\; \Delta'x\rightarrow E\}\; .\]
We put
\[\Phi_{z}(f): \delta'x\mapsto g\cdot f(g^{-1}\delta'x) \; .\]
Observe that this is the place where we actually use that $E$ is a representation of the group $\mathfrak{G}$.
In any case, together with the cohomomorphism $\Phi$, we obtain the induced homomorphism $\varphi^{\ast}$.
This completes the definition of the Hecke operator $T(g)$.

\medskip
The same construction works for cohomology with compact supports. Hence we get also Hecke operators
\[ T(g):H^{\bullet}_{c}(G\setminus X;\mathcal{E})\rightarrow H^{\bullet}_{c}(G\setminus X;\mathcal{E})\; ,\]
\[ T(g):H^{\bullet}_{\textbf{!}}(G\setminus X;\mathcal{E})\rightarrow H^{\bullet}_{\textbf{!}}(G\setminus X;\mathcal{E})\; ,\]
The following properties are easy to verify:
\begin{fact}
 On $H^{0}(G\setminus X;\mathcal{E})=E^G$, the operator $T(g)$ is given by
\[T(g) \cdot e=\sum_{i=1}^{m}\gamma_{i}g\cdot e\]
where the $\gamma_{i}$ are as above. 
\end{fact}

\begin{fact}
 If $\gamma ,\gamma'$ are in $G$, we have
\[T(\gamma g \gamma')=T(g)\; .\]
\end{fact}
The transfer $\tau$ has the following well-known properties:
\begin{itemize}
 \item $\tau \circ \pi'^{\ast}$ is multiplication by the number $[G:\Delta']$.
 \item $\tau:H^{0}(\Delta'\setminus X;k)\rightarrow H^{0}(G\setminus X;k)$ is the homomorphism $[G:\Delta']:k\rightarrow k$.
 \item The following diagram is commutative (use 2.1):
\[\xymatrix{ H^{n}_{c}(\Delta'\setminus X;k) \ar[rr]^{\displaystyle \tau}\ar[rd]_{\int} && H^{n}_{c}(G\setminus X;k)\ar[ld]^{\int} \\
& k  } \] 
 \item $\tau (b \cup \pi'^{\ast}(c))=\tau (b) \cup c\quad $  (cf. [B2], p.405)
 \end{itemize}

From these formulas, we conclude that Hecke operators behave well with respect to the pairing
\[B:H^{m}_{c}(Y;\mathcal{E})\otimes H^{n-m}(Y;\mathcal{E}^{\ast}) \rightarrow k\, .\]
from section 2: We have 
\begin{equation}
 B(T(g)\cdot u \otimes v)= B(u \otimes T(g^{-1}) \cdot v)\; .
\end{equation}
So far in this section, we haven't used any restriction on the characteristics of $k$; it is not even necessary
that $k$ is a field. Now we see:

\begin{theorem}
 Under the assumptions of the duality theorem 2.11, the pair $(H^{m}_{\textbf{!}}(Y;\mathcal{E}),T(g))$ consisting 
of a vector space and an endomorphism is dual to the pair $(H^{n-m}_{\textbf{!}}(Y;\mathcal{E}^{\ast}),T(g^{-1}))$.
\end{theorem}

\end{document}